\documentclass[12pt]{amsart}
\usepackage{amsfonts}
\usepackage{amssymb}
\usepackage{amsmath}
\input xy
\xyoption {all}

\usepackage{geometry}
\usepackage{amsthm}
\usepackage{amssymb}
\usepackage{latexsym}
\usepackage[dvips]{graphics}

\newtheorem{thm}{Theorem}[section] 
\newtheorem{cor}[thm]{Corollary}
\newtheorem{defn}[thm]{Definition}
\newtheorem{example}[thm]{Example}

\newtheorem{lemma}[thm]{Lemma}
\newtheorem{prop}[thm]{Proposition}
\newtheorem{remark}[thm]{Remark}

\newtheorem*{ack}{Acknowledgements}

\numberwithin{equation}{section}

\newcommand{\DC}{{\mathcal D}}

\newcommand{\OO}{{\mathcal O}}
\newcommand{\HC}{{\mathcal H}}
\newcommand{\VV}{{\mathcal V}}

\newcommand{\Z}{\mathbb{Z}}

\newcommand{\Q}{\mathbb{Q}}

\newcommand{\C}{\mathbb{C}}

\newcommand{\PP}{\mathbb{P}}

\newcommand{\VB}{\mathbb{V}}

\begin{document}
\date{\today}
\title{ Hodge genera of algebraic varieties, I.}

\author[Sylvain E. Cappell ]{Sylvain E. Cappell}
\address{S. E. Cappell: Courant Institute, New York University, New York, NY-10012}
\email {cappell@cims.nyu.edu}

\author[Laurentiu G. Maxim ]{Laurentiu G. Maxim}
\address{L. G. Maxim: Courant Institute, New York University, New York, NY-10012}
\email {maxim@cims.nyu.edu}

\author[Julius L. Shaneson ]{Julius L. Shaneson}
\address{J. L. Shaneson: Department of Mathematics, University of Pennsylvania, Philadelphia, PA-19104}
\email {shaneson@sas.upenn.edu}


\thanks{The first and third authors partially supported by grants from NSF and DARPA.
The second author partially supported by a grant from the NYU
Research Challenge Fund.}

\begin{abstract}
The aim of this paper is to study the behavior of Hodge-theoretic
(intersection homology) genera and their associated characteristic
classes under proper morphisms of complex algebraic varieties. We
obtain formulae that relate (parametrized families of) global
invariants of a complex algebraic variety $X$ to such invariants of
singularities of proper algebraic maps defined on $X$. Such formulae
severely constrain, both topologically and analytically, the
singularities of complex maps, even between smooth varieties.
Similar results were announced by the first and third author in
\cite{CS1, S}.
\end{abstract}

\maketitle

\tableofcontents

\section{Introduction}

This paper and its sequels study the behavior of Hodge-theoretic
(intersection homology) genera and their associated characteristic
classes under proper morphisms of complex algebraic varieties. The
formulae obtained in the present paper relate global invariants to
singularities of complex algebraic maps. They thus  shed some light
on the mysterious formulae announced some years ago by the first and
third author in \cite{CS1, S}.

These formulae can be viewed as, on the one hand, yielding powerful
methods of inductively calculating (even parametrized families of)
characteristic classes of algebraic varieties (e.g., by applying
them to resolutions of singularities). On the other hand, they can
be viewed as yielding powerful topological and analytic constraints
on the singularities of any proper algebraic morphism (e.g., even
between smooth varieties), expressed in terms of (even parametrized
families of) their characteristic classes. (Both these perspectives
will be more fully developed following our subsequent studies of the
contributions of monodromy, e.g., see our forthcoming paper
\cite{CLMS}.) Among these severe parameterized constraints on
singularities of maps obtained here in complex settings, only at one
special value do these formulae have full analogues for non-complex
maps (at $y=1$, where they yield topological constraints on the
signature and associated $L$-polynomials of Pontrjagin classes); for
that constraint on topological maps, see Cappell-Shaneson \cite{CS2}
(which employed very different methods) and a comparison in Remark
\ref{one} below.

The main instrument used in proving our results is the functorial
calculus of Grothendieck groups of algebraic mixed Hodge modules.
Originally, some of these results were proven by using Hodge
theoretic aspects of a deep theorem of Bernstein, Beilinson, Deligne
and Gabber, namely the decomposition theorem for the push-forward of
an intersection cohomology complex under a proper algebraic morphism
\cite{BBD,Cat, Cat2}; for details on this approach see \cite{CMS2}.
The functorial approach employed here was suggested to us by the
referee. However, the core calculations used in proving the results
in this paper are modeled on our original considerations based on
BBDG. We assume the reader's familiarity with the language of
sheaves and derived categories, as well as that of intersection
homology, perverse sheaves and Deligne's mixed Hodge structures. But
we do not assume any knowledge about Saito's mixed Hodge modules,
except maybe the now classical notion of an admissible variation of
mixed Hodge structures.


We now briefly outline the content of each section and summarize our
main results. The paper is divided in three main sections: in
Section \S\ref{genera} we discuss genera of complex algebraic
varieties, in \S\ref{Grot} we describe the functorial calculus of
Grothendieck groups of algebraic mixed Hodge modules, while Section
\S\ref{four} deals with characteristic classes yielding the Hodge
genera considered in \S\ref{genera}.

In  \S\ref{fam}, we first recall the definition of Hirzebruch's
$\chi_y$-genus  for a smooth complex projective variety \cite{H},
and explain some of its Hodge theoretic extensions to genera of
possibly singular and/or non-compact complex algebraic varieties,
defined by means of Deligne's mixed Hodge structures. Then we study
the behavior of these genera under proper algebraic maps $f:X \to Y$
which are locally trivial topological fibrations over a compact,
connected, smooth base $Y$. If $F$ denotes the general fiber of $f$,
then, under the assumption of trivial monodromy, $\chi_y$ is
multiplicative, i.e., $\chi_y(X)=\chi_y(Y) \cdot \chi_y(F)$.

In  \S\ref{two}, we consider arbitrary proper algebraic morphisms
$f:X \to Y$ of complex algebraic varieties, and discuss
generalizations of the above multiplicativity property to this
general setting. By taking advantage of the mixed Hodge structure on
the intersection cohomology groups of a possibly singular complex
algebraic variety $X$ \cite{Cat, Sa, Sa1, MHM}, we define
intersection homology genera, $I\chi_y(X)$, that encode the
intersection homology Hodge numbers, and provide yet another
extension of Hirzebruch's genera to the singular case. For example,
$I\chi_{-1}(X)$ is the intersection homology Euler characteristic of
$X$, and if $X$ is projective then $I\chi_1(X)=\sigma(X)$ is the
Goresky-MacPherson signature (\cite{GM1}) of the intersection form
on the middle-dimensional intersection homology of $X$. $I\chi_0(X)$
can be regarded as an extension to singular varieties of the
arithmetic genus. The main results of this section are Theorems
\ref{formula1} and \ref{formula2}, in which the (intersection
homology) genus $\chi_y(X)$ (and resp. $I\chi_y(X)$) of the domain
is expressed in terms of the singularities of the map. More
precisely, we first fix an algebraic stratification of the proper
morphism $f$, that is, we choose complex algebraic Whitney
stratifications of $X$ and $Y$ so that $f$ becomes a stratified
submersion. In particular, the strata satisfy the frontier
condition: $V \cap {\bar W} \neq \emptyset$ implies $V \subset {\bar
W}$. Then the finite set $\VV:=\{V\}$ of all strata of $Y$ is
partially ordered by ``$V \leq W$ if and only if $V \subset {\bar
W}$". If, moreover, we assume that $Y$ is irreducible, then there is
exactly one top-dimensional stratum $S$ of $Y$, with ${\rm dim} S=
{\rm dim} Y$, and $S$ is Zariski-open and dense in $Y$, with $V \leq
S$ for all $V \in \VV$. Let $F$ denote the (generic) fiber of $f$
above $S$. Then we have:
\begin{thm}\label{gen} Let $f :X \to Y$ be a proper algebraic map
of complex algebraic varieties, with $Y$ irreducible. Let $\VV$ be
the set of components of strata of $Y$ in an algebraic
stratification of $f$, and assume $\pi_1(V)=0$ for all $V \in \VV$.
For each $V \in \VV$ with ${\rm dim} (V)< {\rm dim} (Y)$, define
inductively
$$\widehat{I\chi}_y(\bar V)= I\chi_y(\bar V)- \sum_{W < V}
\widehat{I\chi}_y(\bar W) \cdot I\chi_y(c^{\circ}L_{W,V}),$$ where
the sum is over all $W \in \VV$ with $\bar{W} \subset \bar{V}
\setminus V$, and $c^{\circ}L_{W,V}$ denotes the open cone on the
link of  $W$ in $\bar{V}$. Then:
\begin{enumerate}
\item \begin{equation}\label{E1} \chi_y(X)=I\chi_y(Y) \cdot \chi_y(F) + \sum_{V < S} \widehat{I\chi}_y(\bar V) \cdot \left(
\chi_y(F_V) - \chi_y(F) \cdot I\chi_y(c^{\circ}L_{V,Y})
\right),\end{equation} where $F_V$ is the fiber of $f$ above the
stratum $V$.
\item Assume moreover that $X$ is pure dimensional. Then
\begin{equation}\label{E2} I\chi_y(X)=I\chi_y(Y) \cdot I\chi_y(F) + \sum_{V < S}
\widehat{I\chi}_y(\bar V) \cdot \left( I\chi_y(f^{-1}(c^\circ
L_{V,Y})) - I\chi_y(F) \cdot I\chi_y(c^{\circ}L_{V,Y})
\right).\end{equation}
\end{enumerate}
\end{thm}
Let us explain in more detail the $I\chi_y$-terms appearing in the
above formulae. Since $f$ is stratified, the restrictions
$R^jf_*(\Q_X)|_W$ and $R^jf_*(IC'_X)|_W$, for all $j \in \Z$ and $W
\in \VV$, are admissible variations of mixed Hodge structures, where
$IC'_X := IC_X[-{\rm dim}(X)]$ is the shifted intersection
cohomology complex of X. In fact, the first is the classical example
of a ``geometric variation of mixed Hodge structures", whereas the
assertion for the second follows from Saito's theory of algebraic
mixed Hodge modules. Similarly, by Saito's theory, the cohomology
sheaves $\HC^j(IC'_{\bar V})|_W$, for $j \in \Z$ and $W, V \in
\VV$, underly admissible variation of mixed Hodge structures. And
under the assumption that the fundamental group of $W$ is trivial,
these are trivial variations so that the isomorphism classes of the
mixed Hodge structures
\begin{enumerate}
\item $R^jf_*(\Q_X)_w \simeq H^j(\{f=w\};\Q)$,
\item $R^jf_*(IC'_X)_w \simeq IH^j(\{f=w\};\Q)$ for ${\rm
dim}(W)={\rm dim} (Y)$, and \\
$R^jf_*(IC'_X)_w \simeq IH^j(f^{-1}(c^\circ L_{W,Y});\Q)$ for ${\rm
dim}(W) < {\rm dim} (Y)$,
\item $\HC^j(IC'_{\bar V})_w \simeq IH^j(c^\circ L_{W,V};\Q),$
\end{enumerate}
do not depend on the choice of the point $w \in W$. The same is true
for the corresponding $I\chi_y$-genera, even without the assumption
$\pi_1(W) = 0$, because these genera depend only the corresponding
Hodge numbers, which are constant in a variation of mixed Hodge
structures. Also note that the inductive definition of
$\widehat{I\chi}_y(\bar V)$ only depends on the stratified space
$\bar V$ with its induced stratification. Finally, while the
existence of the mixed Hodge structure on the global cohomology
follows already from Deligne's classical work, and the Hodge
structure on the global intersection cohomology of compact algebraic
varieties can be obtained by others methods as in \cite{Cat, Cat2},
the mixed Hodge structures on the intersection cohomology of the
open links $c^\circ L_{V,Y}$ (or of their inverse images) use the
stalk description above, and essentially depend on Saito's theory of
mixed Hodge modules.

In Section \S\ref{Grot}, we explain in detail the functorial
calculus of Grothendieck groups of algebraic mixed Hodge modules,
and prove in this language the main technical result of this paper.
The proof of this result is exactly the same as that of the
corresponding statement in the framework of Grothendieck groups of
constructible sheaves used in our paper \cite{CMS}. This is indeed
the case because, by Saito's work, the calculus of Grothendieck
groups of constructible sheaves completely lifts to the context of
Grothendieck groups of algebraic mixed Hodge modules. Theorem
\ref{gen} above, as well as the formulae for characteristic classes
discussed below, are direct consequence of this main theorem.

In Section \S \ref{four}, we outline the construction of a natural
transformation, $MHT_y$, which, when evaluated at the intersection
cohomology complex $IC'_X$ of a variety $X$, yields a twisted
homology class $IT_y(X)$, whose associated genus for $X$ compact is
$I\chi_y(X)$. The definition uses Saito's theory of mixed Hodge
modules, and is based on ideas of a recent paper of Brasselet,
Sch\"urmann and Yokura \cite{BSY}. If $X$ is a non-singular complex
algebraic variety, then $IT_y(X)$ is the Poincar\'e dual of the
modified Todd class $T_y^*(TX)$ that appears in the generalized
Hirzebruch-Riemann-Roch theorem. For a proper algebraic map $f:X \to
Y$ with $X$ pure-dimensional and $Y$ irreducible, we prove a formula
for the push-forward of the characteristic class $IT_y(X)$ in terms
of characteristic classes of strata of $f$. The main result of this
section can be stated as follows:
\begin{thm}\label{class}
With the notations and assumptions from the above theorem, for each
$V \in \VV$ define inductively
$$\widehat{IT}_y(\bar V)= IT_y(\bar V)- \sum_{W < V} \widehat{IT}_y(\bar W) \cdot I\chi_y(c^{\circ}L_{W,V}),$$
where all homology characteristic classes are regarded in the
Borel-Moore homology of the ambient space $Y$ (with coefficients in
$\Q[y,y^{-1},(1+y)^{-1}]$). Then:
\begin{equation}\label{E3} f_*IT_y(X)=IT_y(Y) \cdot I\chi_y(F) + \sum_{V < S}
\widehat{IT}_y(\bar V) \cdot \left( I\chi_y(f^{-1}(c^\circ L_{V,Y}))
- I\chi_y(F) \cdot I\chi_y(c^{\circ}L_{V,Y}) \right),\end{equation}
where $L_{V,Y}$ is the link of $V$ in $Y$.
\end{thm}

Without the trivial monodromy assumption, the terms in the formulae
of Theorems \ref{gen} and \ref{class} must be written in terms of
genera and respectively characteristic classes with coefficients in
local systems (variations of Hodge structures) on open strata.

The paper ends by discussing some immediate consequences of the
push-forward formula of Theorem \ref{class}.

\bigskip

In a future paper, we will consider the behavior under proper
algebraic maps of $\chi_y$-genera that are defined by using the
Hodge-Deligne numbers of (compactly supported) cohomology groups of
a possibly singular algebraic variety, and deal with non-trivial
monodromy considerations (cf. \cite{CLMS}). We point out that
preliminary results on the Euler characteristics $\chi_{-1}$,
$I\chi_{-1}$, and on the homology MacPherson-Chern classes
(\cite{M}) $T_{-1}=c_* \otimes \Q$, $IT_{-1}=Ic_* \otimes \Q$ of
complex algebraic (resp. compact complex analytic) varieties, have
been already obtained by the authors in \cite{CMS}.

\begin{ack} We are grateful to the anonymous referee for many
valuable comments, and for suggesting the functorial approach
employed in this paper. We also thank J\"org Sch\"urmann and Mark
Andrea de Cataldo for helpful comments and discussions.
\end{ack}

\section{Genera}\label{genera} In this section, we define
Hodge-theoretic genera of complex algebraic varieties, and study
their behavior under proper algebraic morphisms.

\subsection{Families over a smooth base.}\label{fam}

\begin{defn}\label{D1}\rm
For a smooth projective variety $X$, we define its Hirzebruch
$\chi_y$-genus (\cite{H}) by the formula: \begin{equation}\label{H}
\chi_y(X)=\sum_p \left( \sum_{q} (-1)^{q} h^{p,q}(X) \right) y^p,
\end{equation} where $h^{p,q}(X)={\rm dim}_{\C} H^q(X;\Omega^p_X)$ are the Hodge numbers of $X$.
Note that $\chi_{-1}$ is the usual Euler characteristic, $\chi_0$ is
the arithmetic genus, and $\chi_1$ is the signature of $X$.
\end{defn}

Various extensions of Hirzebruch's genus to the singular and/or
non-compact setting shall be explained below. First note that the
Grothendieck group $K_0(mHs^{(p)})$ of the abelian category of
(graded polarizable) rational mixed Hodge structures is a ring with
respect to the tensor product, with unit the pure Hodge structure
$\Q$ of weight zero. We can now make the following
\begin{defn}\label{D2}\rm The $E$-polynomial is the ring homomorphism $$E:
K_0(mHs^{(p)}) \to \Z[u,v,u^{-1},v^{-1}]$$ defined by
\begin{equation} (V,F^{\bullet},W_{\bullet}) \mapsto E(V):=\sum_{p,q} {\rm
dim}_{\C} (gr^p_F gr^W_{p+q}(V \otimes_{\Q} \C)) \cdot u^pv^q.
\end{equation} This is well-defined since the functor $gr^p_F
gr^W_{p+q}(- \otimes_{\Q} \C)$ is  exact on mixed Hodge structures.
Specializing to $(u,v)=(-y,1)$, we get the $\chi_y$-genus
\begin{equation} \chi_y:K_0(mHs^{(p)}) \to \Z[y,y^{-1}]; \
[(V,F^{\bullet},W_{\bullet})] \mapsto \sum_p {\rm dim}_{\C} (gr^p_F
(V \otimes_{\Q} \C))\cdot (-y)^p.\end{equation}  Using Deligne's
mixed Hodge structure on the cohomology (with compact support)
$H^j_{(c)}(X;\Q)$ of a complex algebraic variety $X$ (cf.
\cite{De}), we can define \begin{equation}\label{p14}
[H^*_{(c)}(X;\Q)]:=\sum_j(-1)^j \cdot [H^j_{(c)}(X;\Q)] \in
K_0(mHs^p).\end{equation}  Then by applying one of the homomorphisms
$E$ or $\chi_y$, we obtain $E(X)$, $E_c(X)$, $\chi_y(X)$ and
$\chi_y^c(X)$, so that for $X$ smooth and projective, this
definition of $\chi_y(X)=\chi_y^c(X)$ agrees with the classical
Hirzebruch genus of $X$ from Definition \ref{D1}.
\end{defn}

We first show that if $f:X \to Y$ is a family of compact varieties
(i.e., a locally trivial topological fibration in the complex
topology) over a smooth, connected, compact variety $Y$, then under
certain assumptions on monodromy, $\chi_y$ behaves multiplicatively.
In the setting of algebraic geometry, this fact encodes as special
cases the classical multiplicativity property of the
Euler-Poincar\'{e} characteristic for a locally trivial topological
fibration, and the Chern-Hirzebruch-Serre formula for the signature
of fibre bundles with trivial monodromy.

\begin{prop}\label{smooth} Let $f:X \to Y$ be a proper algebraic map
of complex algebraic varieties, with $Y$ compact, smooth and
connected. Suppose that all direct image sheaves $R^jf_*\Q_X$, $j
\in \Z$, are locally constant, e.g., $f$ is a locally trivial
topological fibration. Assume $\pi_1(Y)$ acts trivially on the
cohomology of the general fiber $F$ of $f$ (e.g., $\pi_1(Y)=0$),
i.e., all these local systems $R^jf_*\Q_X$ ($j \in \Z$) are
constant. Then:
\begin{equation}\label{E10} \chi_y(X)=\chi_y(Y) \cdot \chi_y(F).\end{equation}
\end{prop}

\begin{proof} The local systems $R^jf_*\Q_X$ ($j \in \Z$) define geometric
variations of mixed Hodge structures, thus admissible variations in
the sense of Steenbrink-Zucker and Kashiwara (cf. \cite{K,SZ}).
Since $Y$ is a smooth, compact algebraic variety, the cohomology
groups $H^i(Y; R^jf_*\Q_X)$ get an induced (polarizable) mixed Hodge
structure. So, for each $j \in \Z$ we can define
$$[H^*(Y;R^jf_*\Q_X)]:=\sum_i (-1)^i \cdot [H^i(Y;R^jf_*\Q_X)] \in
K_0(mHs^p).$$ Moreover, the following key equality holds
\begin{equation}\label{E11}
[H^*(X;\Q)]=\sum_j (-1)^j [H^*(Y;R^jf_*\Q_X)] \in K_0(mHs^p).
\end{equation}
The proof of this equality will be given in Proposition
\ref{gsmooth} of Section \S\ref{calc}, in terms of mixed Hodge
modules. Finally, if the local system $R^jf_*\Q_X)$ is constant, we
have an isomorphism of mixed Hodge structures
$$H^i(Y;R^jf_*\Q_X) \cong H^i(Y;\Q) \otimes H^j(F;\Q).$$
Altogether, we get the equality
\begin{equation}\label{p13}[H^*(X;\Q)]=[H^*(Y;\Q)] \cdot [H^*(F;\Q)]
\in  K_0(mHs^p).\end{equation} We obtain the claimed
multiplicativity by applying  the $\chi_y$-genus homomorphism to the
identity in formula (\ref{p13}). Note that by applying the
$E$-polynomial to (\ref{p13}), we obtain a similar multiplicativity
property for the $E$-polynomials.

\end{proof}

\subsection{Proper maps of complex algebraic varieties}\label{two}
Let $f :X \to Y$ be a proper map of complex algebraic varieties.
Such a map can be stratified with subvarieties as strata. In
particular, there is a filtration of $Y$ by closed subvarieties,
underlying a Whitney stratification $\VV$, so that the restriction
of $f$ to the preimage of any component of a stratum in $Y$ is a
locally trivial map of Whitney stratified spaces, i.e., $f$ becomes
a stratified submersion.

In this paper all intersection cohomology complexes are those
associated to the middle perversity. By convention, the restriction
of the intersection cohomology complex $IC_X$ to the dense open
stratum of $X$ is the constant sheaf shifted by the complex
dimension of $X$. If $X$ is a complex algebraic variety of pure
dimension $n$, the intersection cohomology groups (with compact
support) are defined by the rule
$IH_{(c)}^k(X;\Q):=H_{(c)}^{k-n}(X;IC_X)$.

Another possible extension of Hirzebruch's $\chi_y$-genus to the
singular setting is obtained by using intersection homology theory
(\cite{GM1,GM2} as follows:
\begin{defn}\rm For a pure dimensional complex algebraic variety $X$ we let
$$IC'_X:=IC_X[-{\rm dim}(X)]$$ be the shifted intersection cohomology
complex. Then, by Saito's theory of mixed Hodge modules, the group
$IH_{(c)}^k(X;\Q):=H_{(c)}^{k}(X;IC'_X)$ gets a (graded polarizable)
mixed Hodge structure, which is pure of weight $k$ if $X$ is
compact. So we can define \begin{equation}[IH^*_{(c)}(X;\Q)]:=\sum_j
(-1)^j [IH^j_{(c)}(X;\Q)] \in K_0(mHs^p),\end{equation} and
polynomials
\begin{equation}IE_{(c)}(X):=E([IH^*_{(c)}(X;\Q)])\end{equation} and
\begin{equation}\label{p8}I\chi_y^{(c)}:=\chi_y([IH^*_{(c)}(X;\Q)]).\end{equation}
\end{defn}

As a natural extension of the multiplicativity property for
$\chi_y$-genera of families over a smooth base described in
Proposition \ref{smooth}, we aim to find the deviation from
multiplicativity of the $\chi_y$- and resp. $I\chi_y$-genus in the
more general setting of an arbitrary proper algebraic map. The
formulae proved in Theorem \ref{formula1} and resp. Theorem
\ref{formula2} below include correction terms corresponding to
genera of strata and of their normal slices, and to those of fibers
of $f$ above each stratum in $Y$. The first main result of this
section concerns the $\chi_y$-genus:
\begin{thm}\label{formula1}
Let $f :X \to Y$ be a proper algebraic map of complex algebraic
varieties, with $Y$ irreducible. Let $\VV$ be the set of components
of strata of $Y$ in an algebraic stratification of $f$, and assume
$\pi_1(V)=0$ for all $V \in \VV$.\footnote{Contributions of
non-trivial monodromy to such formulae will be subject of further
studies, e.g., see our forthcoming paper \cite{CLMS}; see also
\cite{BCS, Ba} for some results on monodromy contributions for
signatures and related characteristic classes in topological
settings.} For each $V \in \VV$ with ${\rm dim} (V)< {\rm dim} (Y)$,
define inductively
$$\widehat{I\chi}_y(\bar V)= I\chi_y(\bar V)- \sum_{W < V}
\widehat{I\chi}_y(\bar W) \cdot I\chi_y(c^{\circ}L_{W,V}),$$ where
the sum is over all $W \in \VV$ with $\bar{W} \subset \bar{V}
\setminus V$, and $c^{\circ}L_{W,V}$ denotes the open cone on the
link of  $W$ in $\bar{V}$. Then: \begin{equation}\label{E20}
\chi_y(X)=I\chi_y(Y) \cdot \chi_y(F) + \sum_{V < S}
\widehat{I\chi}_y(\bar V) \cdot \left( \chi_y(F_V) - \chi_y(F) \cdot
I\chi_y(c^{\circ}L_{V,Y}) \right),\end{equation} where $F$ is the
(generic) fiber over the top-dimensional stratum $S$, and $F_V$ is
the fiber of $f$ above the stratum $V \in \VV \setminus \{S\}$.
\end{thm}

The proof is a direct consequence of the calculus of Grothendieck
groups of algebraic mixed Hodge modules, and will be given in
Section \S\ref{MAIN}, Proposition \ref{gsmp}. We want to point out
that similar formulae also hold for $\chi_y^c(X)$ and resp.
$E_{(c)}(X)$, which in the same way follow from the results of
Section \ref{MAIN}.

In the special case when $f$ is the identity map, the equation
(\ref{E20}) measures the difference between the $\chi_y$- and the
$I\chi_y$-genus of a complex algebraic variety (with no monodromy
restrictions in the case of Euler characteristics, that is, for
$y=-1$; see \cite{CMS}, Cor. 3.5):
\begin{cor} Let $Y$ be an irreducible complex algebraic variety, and assume
there is a Whitney stratification $\VV$ of $Y$ with all strata
simply connected. Then, in the notations of the above theorem, we
obtain
\begin{equation}
\chi_y(Y)=I\chi_y(Y) + \sum_{V < S} \widehat{I\chi}_y(\bar V) \cdot
\left( 1 -  I\chi_y(c^{\circ}L_{V,Y}) \right).
\end{equation}
\end{cor}

\begin{example}\label{bup}\rm \emph{Smooth blow-up}\newline
Let $Y$ be a smooth compact $n$-dimensional variety and $Z \subset
Y$ a submanifold of pure codimension $r+1$. Let $X$ be the blow-up
of $Y$ along $Z$, and $f:X \to Y$ be the blow-up map. Then $X$ is a
$n$-dimensional smooth variety, and $f$ is an isomorphism over $Y
\setminus Z$ and a projective bundle (Zariski locally trivial with
fibre $\C\PP^r$) over $Z$, corresponding to the projectivization of
the normal bundle of $Z$ in $Y$ of rank $r+1$. As we later explain
in Example \ref{MainEx}, the formula (\ref{E20}) of Theorem
\ref{formula1} is also true in this context (without assuming that
all strata $V$ are simply-connected) and it reduces to a more
familiar one (\cite{BSY}, Example 3.3):
\begin{equation}\label{BSY}
\chi_y(X)=\chi_y(Y)+\chi_y(Z) \cdot \left( -y+\cdots +(-y)^r
\right).\end{equation} In fact formula (\ref{BSY}) can be easily
obtained just by using the (additivity and multiplicativity)
properties of the $\chi_y^c$-genera of complex algebraic varieties
(cf. \cite{D}), and it holds if one considers $X$ to be the blow-up
of a complete variety $Y$ along a \emph{regularly} embedded
subvariety $Z$ of pure codimension $r+1$. By the "Weak Factorization
Theorem" \cite{AKMW}, any birational map $h:S \to T$ between
complete non-singular complex algebraic varieties can be decomposed
as a finite sequence of projections from smooth spaces lying over
$T$, which are obtained by blowing up or blowing down along smooth
centers. Then  (\ref{BSY}) yields the birational invariance of the
arithmetic genus $\chi_0$ of non-singular projective varieties (see
the discussion in \cite{BSY}, Example 3.3).

\end{example}

\begin{remark}\rm Formula (\ref{E20}) yields calculations of classical topological
and algebraic invariants of the complex algebraic variety $X$, e.g.
Euler characteristic, and if $X$ is smooth and projective, signature
and arithmetic genus, in terms of singularities of proper algebraic
maps defined on $X$.

\end{remark}

The second main result of this section asserts that the
$I\chi_y$-genus defined in (\ref{p8}) satisfies the so-called
``stratified multiplicative property" (compare \cite{CS1,S}). More
precisely, we have the following:

\begin{thm}\label{formula2}
Let $f :X \to Y$ be a proper algebraic map of complex algebraic
varieties, with $X$ pure dimensional and $Y$ irreducible. Let $\VV$
be the set of components of strata of $Y$ in an algebraic
stratification of $f$, and assume $\pi_1(V)=0$ for all $V \in \VV$.
Then in the notations of Theorem \ref{formula1}, we have
\begin{equation}\label{E21} I\chi_y(X)=I\chi_y(Y) \cdot I\chi_y(F) + \sum_{V < S}
\widehat{I\chi}_y(\bar V) \cdot \left( I\chi_y(f^{-1}(c^\circ
L_{V,Y})) - I\chi_y(F) \cdot I\chi_y(c^{\circ}L_{V,Y})
\right).\end{equation}
\end{thm}

Again, the proof will follow from the considerations of Section
\ref{MAIN}, Prop. \ref{gsmp}, which will also imply that similar
formulae hold for $I\chi_y^c(X)$ and $IE_{(c)}(X)$. Let us only
explain here the identification of stalks at a point $w$ in a
stratum $W$ of $Y$, that yield the $I\chi_y$-terms in formula
(\ref{E21}). Let $i_w:\{w\} \hookrightarrow Y$ be the inclusion map.
Then

\begin{lemma}\label{id} We have the following isomorphisms
\begin{enumerate}
\item $R^jf_*(IC'_X)_w = \HC^j(i_w^*Rf_*(IC'_X)) \simeq IH^j(\{f=w\};\Q)$ for ${\rm
dim}(W)={\rm dim} (Y)$, and \\
$R^jf_*(IC'_X)_w = \HC^j(i_w^*Rf_*(IC'_X)) \simeq
IH^j(f^{-1}(c^\circ L_{W,Y});\Q)$ for ${\rm dim}(W) < {\rm dim}
(Y),$
\item $\HC^j(IC'_{\bar V})_w = \HC^j(i_w^*IC'_{\bar V}) \simeq IH^j(c^\circ L_{W,V};\Q),$
\end{enumerate}
which endow the intersection cohomology groups on the right hand
side of the above identities with canonical mixed Hodge structures.
\end{lemma}
\begin{proof}
If ${\rm dim}(W)={\rm dim} (Y)$, then $\{f=w\}$ is the generic fiber
$F$ of $f$, thus it is locally normally non-singular embedded in
$X$.  It follows from [\cite{GM2}, \S 5.4.1] that we have a
quasi-isomorphism:
$$IC'_X |_{F} \simeq IC'_{F}.$$
Then by proper base change, we obtain that $\HC^j(i_w^*Rf_*(IC'_X))
\simeq IH^j(\{f=w\};\Q)$.

Assume now that ${\rm dim}(W) < {\rm dim} (Y)$. Let $N$ be a normal
slice to $W$ at $w$ in local analytic coordinates $(Y,w)
\hookrightarrow (\C^n,w)$, that is, a germ of a complex manifold
$(N,w) \hookrightarrow (\C^n,w)$, intersecting $W$ transversally
only at $w$, and with ${\rm dim} (W) + {\rm dim} (N)=n$. Recall that
the link $L_{W,Y}$ of the stratum $W$ in $Y$ is defined as
$$L_{W,Y}:=Y \cap N \cap
\partial B_r(w),$$
where $B_r(w)$ is an open ball of (very small) radius $r$ around
$w$. Moreover, $Y \cap N \cap B_r(w)$ is isomorphic (in a stratified
sense) to the open cone $c^{\circ}L_{W,Y}$ on the link (\cite{B}, p.
44).  By factoring the inclusion map $i_w$ as the composition:
$$\{w\} \overset{\phi}{\hookrightarrow} Y \cap N
\overset{\psi}{\hookrightarrow} Y$$ we can now write:
{\allowdisplaybreaks
\begin{eqnarray*}
\HC^j(i_w^*Rf_*IC'_X)
&\cong& H^j(w,\phi^*\psi^*Rf_*IC'_X)\\
&\cong& H^j(\psi^*Rf_*IC'_X)_{w}\\
&\cong& H^j(c^{\circ}L_{W,Y},Rf_*IC'_X)\\
&\cong& H^j(f^{-1}(c^{\circ}L_{W,Y}),IC'_X)\\
&\overset{(1)}{\cong}& H^j(f^{-1}(c^{\circ}L_{W,Y}),IC'_{f^{-1}(c^{\circ}L_{W,Y})})\\
&\cong& IH^{j}(f^{-1}(c^{\circ}L_{W,Y});\Q)
\end{eqnarray*}
} where in $(1)$ we used the fact that the inverse image of a normal
slice to a stratum of $Y$ in a stratification of $f$ is (locally)
normally non-singular embedded in $X$ (this fact is a consequence of
first isotopy lemma).

The isomorphism $\HC^j(i_w^*IC'_{\bar V}) \simeq IH^j(c^\circ
L_{W,V};\Q)$ follows, for example, from [\cite{B}, p.30, Prop.4.2].

The statement about the existence of the canonical mixed Hodge
structures is a consequence of Saito's theory of mixed Hodge modules
(see Section \ref{smhm}). Indeed, since the complexes $Rf_*IC'_X$
and resp. $IC'_{\bar V}$ underlie complexes of mixed Hodge modules
(cf. \S\ref{Grot}), their pullbacks over the point $w$ become
complexes of mixed Hodge structures, so their cohomologies are
(graded polarizable) rational mixed Hodge structures.

\end{proof}

\begin{remark}\label{one}\rm For a projective algebraic variety $X$ of pure dimension $n$,
the value at $y=1$ of the intersection homology genus $I\chi_y(X)$
is the Goresky-MacPherson signature $\sigma(X)$ of the intersection
form in the middle-dimensional intersection cohomology $IH^n(X; \Q)$
with middle-perversity (\cite{GM1}). Therefore, under the trivial
monodromy assumption, formula (\ref{E21}) calculates the signature
of the domain of a proper map $f$ in terms of singularities of the
map. In \cite{CS2}, a different formula was given for the behavior
of the signature (and associated $L$-class) under any stratified
map. Those topological results were obtained by a very different
sheaf theoretic method, i.e., introducing a notion of cobordism of
self-dual sheaves and showing sheaf decompositions up to such
cobordism. By comparing the two formulae in the case of a proper map
of algebraic varieties, one obtains interesting Hodge theoretic
interpretations of the normal data encoded in the topological
formula for signature \cite{CS2}. We exemplify this relation in a
simple situation, namely that of blowing up a point: Let $X$ be
obtained from $Y$ by blowing up a point $y$. Let $L$ be the link of
$y$ in $Y$. Formula (\ref{E21}) becomes in this case:
$$\sigma(X)=\sigma(Y)+I\chi_1(f^{-1}(c^{\circ}L))-I\chi_1(c^{\circ}L).$$
On the other hand, the topological formula for signature in \cite{CS2} yields
$$\sigma(X)=\sigma(Y)+\sigma(E_y),$$ where $E_y=f^{-1}(N)/f^{-1}(L)$ is the topological
completion of $f^{-1}(\text{int}N)$, for $N$ a piecewise linear
neighborhood of $y$ in $Y$ with $\partial N=L$. By comparing the two
formulae above, we obtain a Hodge theoretic interpretation for the
signature of the topological completion $E_y$, namely:
$$\sigma(E_y)=I\chi_1(f^{-1}(c^{\circ}L))-I\chi_1(c^{\circ}L).$$

\end{remark}

\section{Grothendieck groups of algebraic mixed Hodge
modules}\label{Grot} In this section we explain the functorial
calculus of Grothendieck groups of algebraic mixed Hodge modules,
and prove in this language the main technical result of this paper.
All other results, in particular those stated in Section
\ref{genera}, are direct consequences of this main theorem.

\subsection{Mixed Hodge Modules}\label{smhm}  For ease of reading,
we begin by introducing a few notions that will be used throughout
this paper. A quick introduction to Saito's theory can be also found
in \cite{BSY,PS,Sa}. Generic references are \cite{Sa0, Sa1}, but see
also \cite{MHM}.

Let $X$ be an $n$-dimensional complex algebraic variety. To such an
$X$ one can associate an abelian category of algebraic mixed Hodge
modules, $MHM(X)$, together with  functorial push-downs $f_!$, $f_*$
on the level of derived categories $D^bMHM(X)$ for any, not
necessarily proper, map. If $f$ is a proper map, then $f_*=f_!$. In
fact, the derived category $D^b_c(X)$ of bounded (algebraically)
constructible complexes of sheaves of $\Q$-vector spaces
\emph{underlies} the theory of mixed Hodge modules, i.e., there is a
forgetful functor
$$rat: D^bMHM(X) \to D^b_c(X)$$ which associates their underlying
$\Q$-complexes to complexes of mixed Hodge modules, so that $$rat (
MHM(X)) \subset Perv(\Q_X),$$ that is to say that $rat \circ H =
{^{p}\HC} \circ rat$, where $H$ stands for the cohomological functor
in $D^bMHM(X)$ and ${^{p}\HC}$ denotes the perverse cohomology. Then
the functors $f_*$, $f_!$, $f^*$, $f^!$, $\otimes$, $\boxtimes$ on
$D^bMHM(X)$ are ``lifts" of the similar functors defined on
$D^b_c(X)$, with $(f^*,f_*)$ and $(f_!,f^!)$ also pairs of adjoint
functors in the context of mixed Hodge modules.

The objects of the category $MHM(X)$ can be roughly described as
follows. If $X$ is \emph{smooth}, then $MHM(X)$ is a full
subcategory of the category of objects $((M,F), K, W)$ such that:
(1)  $(M,F)$ is an algebraic holonomic filtered $\DC$-module $M$ on
$X$, with an exhaustive, bounded from below and increasing ``Hodge"
filtration $F$ by algebraic $\OO_X$-modules; (2)  $K \in Perv(\Q_X)$
is the underlying rational sheaf complex, and there is a
quasi-isomorphism $\alpha: DR(M) \simeq \C \otimes K$ in
$Perv(\C_X)$, where $DR$ is the de Rham functor shifted by the
dimension of $X$; (3)  $W$ is a pair of filtrations on $M$ and $K$
compatible with $\alpha$. For a singular $X$, one works with
suitable local embeddings into manifolds and corresponding filtered
$\DC$-modules with support on $X$. In addition, these objects have
to satisfy a long list of very complicated properties, but the
details of the full construction are not needed here. Instead, we
will only use certain formal properties that will be explained
below. In this notation, the functor $rat$ is defined by $rat((M,F),
K, W)=K$.

It follows from the definition of mixed Hodge modules that every $M
\in MHM(X)$ has a functorial increasing filtration $W$ in $MHM(X)$,
called the \emph{weight filtration} of $M$, so that the functor $M
\to Gr_k^W M$ is exact. We say that $M \in MHM(X)$ is \emph{pure of
weight $k$} if $Gr_i^W M = 0$ for all $i \neq k$. A complex
$M^{\bullet} \in D^bMHM(X)$ is \emph{mixed of weight $\leq k$ (resp.
$\geq k$)} if $Gr_i^W H^jM^{\bullet} = 0$ for all $i
> j+k$ (resp. $i < j+k$), and it is \emph{pure of weight $k$} if
$Gr_i^W H^jM^{\bullet} = 0$ for all $i \neq j+k$. If $f$ is a map of
algebraic varieties, then $f_!$ and $f^*$ preserve weight $\leq k$,
and $f_*$ and $f^!$ preserve weight $\geq k$. If $M^{\bullet} \in
D^bMHM(X)$ is of weight $\leq k$ (resp. $\geq k$), then
$H^jM^{\bullet}$ has weight $\leq j+k$ (resp. $\geq j+k$). In
particular, if $M \in D^bMHM(X)$ is pure and $f:X \to Y$ is proper,
then $f_*M \in D^b(MHM(Y))$ is pure.

If $j: U \hookrightarrow X$ is a Zariski-open dense subset in $X$,
then the intermediate extension $j_{!*}$ (cf. \cite{BBD}) preserves
the weights.

We say that $M \in MHM(X)$ is supported on $S$ if and only if
$rat(M)$ is supported on $S$. Saito showed that the category of
mixed Hodge modules supported on a point, $MHM(pt)$, coincides with
the category $mHs^p$ of (graded) polarizable rational mixed Hodge
structures. Here one has to switch the increasing
$\mathcal{D}$-module filtration $F_*$ of the mixed Hodge module to
the decreasing Hodge filtration of the mixed Hodge structure by
$F^{*}:=F_{-*}$, so that $gr^p_F \simeq  gr^F_{-p}$. In this case,
the functor $rat$ associates to a mixed Hodge structure the
underlying rational vector space. Following \cite{Sa1}, there exists
a unique object $\Q^H \in MHM(point)$ such that $rat(\Q^H)=\Q$ and
$\Q^H$ is of type $(0,0)$. In fact, $\Q^H=((\C,F), \Q, W)$, with
$gr^F_i=0=gr^W_i$ for all $i \neq 0$, and $\alpha : \C \to \C
\otimes \Q$ the obvious isomorphism. For a complex algebraic variety
$X$, we define the complex of mixed Hodge modules
$$\Q_X^H:=k^*\Q^H \in D^bMHM(X)$$ with $rat(\Q^H_X)=\Q_X$, where
$k:X\to pt$ is the constant map to a point. If $X$ is \emph{smooth}
and of dimension $n$, then $\Q_X[n] \in Perv(\Q_X)$, and
$\Q_X^H[n]\in MHM(X)$ is a single mixed Hodge module (in degree $0$)
explicitly described by $\Q_X^H[n]=((\OO_X, F), \Q_X[n], W),$ where
$F$ and $W$ are trivial filtrations so that $gr^F_i=0=gr^W_{i+n}$
for all $ i \neq 0$. So if $X$ is smooth of dimension $n$, then
$\Q_X^H[n]$ is pure of weight $n$. By the stability of the
intermediate extension functor, this shows that if $X$ is an
algebraic variety of pure dimension $n$ and $j:U \hookrightarrow Z$
is the inclusion of a smooth Zariski-open dense subset, then the
intersection cohomology module $IC_X^H:=j_{!*}(\Q_U^H[n])$ is pure
of weight $n$, with underlying perverse sheaf $rat(IC_X^H)=IC_X$.

The pure objects in Saito's theory are the \emph{polarized Hodge
modules} (see \cite{Sa0}, but here we work in the more restricted
algebraic context as defined in \cite{Sa1}). The category
$MH(X,k)^p$ of polarizable Hodge modules on $X$ of weight $k$ is a
semi-simple abelian category, in the sense that every polarizable
Hodge module on $X$ can be written in a unique way as a direct sum
of polarizable Hodge modules with strict support in irreducible
closed subvarieties of $X$. This is the so-called
\emph{decomposition by strict support} of a pure Hodge module. If
$MH_Z(X,k)^p$ denotes the category of pure Hodge modules of weight
$k$ and strict support $Z$, then $MH_Z(X,k)^p$ depends only on $Z$,
and any $M \in MH_Z(X,k)^p$ is \emph{generically} a polarizable
variation of Hodge structures $\VB_U$ on a Zariski open dense subset
$U \subset Z$, with \emph{quasi-unipotent monodromy at infinity}.
Conversely, any such polarizable variation of Hodge structures can
be extended uniquely to a pure Hodge module. In other words, there
is an equivalence of categories: \begin{equation}MH_Z(X,k)^p \simeq
VHS_{gen}(Z, k- dim(Z))^p,\end{equation} where the right-hand side
is the category of polarizable variations of Hodge structures of
weight $k-{\rm dim}(Z)$ defined on non-empty smooth subvarieties of
$Z$, whose local monodromies are quasi-unipotent. Note that, under
this correspondence, if $M$ is a pure Hodge module with strict
support $Z$ then $rat(M)=IC_Z(\VB)$, where $\VB$ is the
corresponding variation of Hodge structures.

If $X$ is smooth of dimension $n$, an object $M \in MHM(X)$ is
called \emph{smooth} if and only if $rat(M)[-n]$ is a local system
on $X$.  Smooth mixed Hodge modules are (up to a shift) admissible
(at infinity) variations of mixed Hodge structures. Conversely, an
admissible variation of mixed Hodge structures $\VB$ on a smooth
variety $X$ of pure dimension $n$ gives rise a smooth mixed Hodge
module (cf. \cite{Sa1}), i.e., to an element $\VB^H[n]\in MHM(X)$
with $rat(\VB^H[n])=\VB[n]$. A pure polarizable variation of weight
$k$ yields an element of $MH(X,k+n)^p$. By the stability by the
intermediate extension functor it follows that if $X$ is an
algebraic variety of pure dimension $n$ and $\mathbb{V}$ is an
admissible variation of (pure) Hodge structures (of weight $k$) on a
Zariski-open dense subset $U \subset X$, then $IC^H_X(\mathbb{V})$
is an algebraic mixed Hodge module (pure of weight $k+n$), so that
$rat(IC_X^H(\VB)|_U)=\VB[n]$.

We conclude this section with a short explanation of the ``rigidity"
property for admissible variations of mixed Hodge structures, as
this will be used later on. Assume $X$ is smooth, connected and of
dimension $n$, with $M \in MHM(X)$ a smooth mixed Hodge module, so
that the local system $\VB:=rat(M)[-n]$ has the property that the
restriction map $r:H^0(X;\VB) \to \VB_x$ is an isomorphism for all
$x \in X$. Then the (admissible) variation of mixed Hodge structures
is a constant variation since $r$ underlies the morphism of mixed
Hodge structures (induced by the adjunction $id \to i_*i^*$):
$$H^0(k_*(M)[-n]) \to H^0(k_*i_*i^*(M)[-n])$$
with $k:X \to pt$ the constant map, and $i :\{x\} \hookrightarrow X$
the inclusion of the point. This implies $$M[-n]=\VB^H\simeq
k^*i^*\VB^H=k^*\VB^H_x \in D^bMHM(X).$$

\subsection{The functorial calculus of Grothendieck
groups.}\label{calc} In this section, we describe the functorial
calculus of Grothendieck groups of algebraic mixed Hodge modules. As
a first application, we indicate the proof of the identity
(\ref{E11}), which was used in Proposition \ref{smooth}.

Let $X$ be a complex algebraic variety. By associating to (the class
of) a complex the alternating sum of (the classes of) its cohomology
objects, we obtain the following identification (e.g. compare
[\cite{KS}, p. 77], [\cite{Sc}, Lemma 3.3.1])
\begin{equation} K_0(D^bMHM(X))=K_0(MHM(X)).
\end{equation}
In particular, if $X$ is a point,
\begin{equation} K_0(D^bMHM(pt))=K_0(mHs^p),
\end{equation}
and the latter is a commutative ring with respect to the tensor
product, with unit $\Q^H_{pt}$. All functors $f_*$, $f_!$, $f^*$,
$f^!$, $\otimes$, $\boxtimes$ induce corresponding functors on
$K_0(MHM(\cdot))$. Moreover, $K_0(MHM(X))$ becomes a
$K_0(MHM(pt))$-module, with the multiplication induced by the exact
exterior product
$$\boxtimes : MHM(X) \times MHM(pt) \to MHM(X \times \{pt\}) \simeq
MHM(X).$$ Also note that $$M \otimes \Q^H_X \simeq M \boxtimes
\Q^H_{pt} \simeq M$$ for all $M \in MHM(X)$. Therefore,
$K_0(MHM(X))$ is a unitary $K_0(MHM(pt))$-module. Finally, the
functors $f_*$, $f_!$, $f^*$, $f^!$ commute with exterior products
(and $f^*$ also commutes with the tensor product $\otimes$), so that
the induced maps at the level of Grothendieck groups
$K_0(MHM(\cdot))$ are $K_0(MHM(pt))$-linear. Moreover, by the
functor $$rat:K_0(MHM(X)) \to K_0(D^b_c(X)) \simeq
K_0(Perv(\Q_X)),$$ these lift the corresponding transformations from
the (topological) level of Grothendieck groups of constructible (or
perverse) sheaves.

\bigskip

As a first application, we can now explain the proof of Proposition
\ref{smooth} in the following more general form.

\begin{prop}\label{gsmooth}
Let $f:X \to Y$ be a proper algebraic map, with $Y$ a smooth,
connected variety, of dimension $n$. Fix $M \in D^bMHM(X)$ so that
all higher direct image sheaves $R^jf_*rat(M)$ ($j \in \Z$) are
locally constant. Then we have
\begin{equation}\label{gE11}
[H^*(X;rat(M)] =\sum_j(-1)^j \cdot [H^*(Y;R^jf_*rat(M))] \in
K_0(mHs^p).
\end{equation}
If, moreover, the local systems $R^jf_*rat(M)$ ($j \in \Z$) are
constant, then we obtain the following multiplicative relation in
$K_0(mHs^p)$ (extending (\ref{p13})):
\begin{equation}\label{gmult}
[H_{(c)}^*(X;rat(M)]=[H_{(c)}^*(Y;\Q)] \cdot [H^*(\{f=y\};rat(M))].
\end{equation}
\end{prop}

\begin{proof} We first need to explain the various mixed Hodge
structures appearing in the above formulae. Let $k:X \to pt$ be the
constant map. Then the cohomology groups
$$H^j(X;rat(M))=rat(H^j(k_*M)) \ \ \text{and} \ \ H^j_{c}
(X;rat(M))=rat(H^j(k_!M))$$ get induced (graded polarizable) mixed
Hodge structures, so that we can define
$$[H^*_{(c)}(X;rat(M))]:=k_{*(!)}[M] \in K_0(mHs^p).$$
By a deep theorem of Saito (cf. \cite{Sa5}), these mixed Hodge
structures for $M=\Q^H_X$ agree with those of Deligne, so that in
this case we get back our old notation from (\ref{p14}).

Since all direct image sheaves $R^jf_*rat(M)$ are locally constant
($j \in \Z$), the usual and perverse cohomology sheaves of
$Rf_*rat(M)$ are the same up to a shift by $n$ (since $Y$ is
smooth), so that for each $j \in \Z$, $H^j(f_*M)$ is a \emph{smooth}
mixed Hodge module with
$$rat(H^{j+n}(f_*M))=R^jf_*rat(M)[n].$$
In particular, the local system $R^jf_*rat(M)$ underlies an
admissible variation of mixed Hodge structures. By pushing down
under the constant map $k':Y \to pt$, we obtain
$$[H^*(Y;R^jf_*rat(M))]=(-1)^n \cdot k'_*([H^{j+n}(f_*M)]) \in
K_0(mHs^p),$$ which, by construction, agrees in the case when
$M=\Q^H_X$ and $Y$ compact with the element of $K_0(mHs^p)$
appearing in the equation (\ref{E11}). By taking the alternating sum
of cohomology sheaves, we also have that
$$f_*[M]=[f_*M]=\sum_j(-1)^j \cdot [H^j(f_*M)] \in K_0(MHM(Y)),$$
and, by functoriality, we get the equality in the formula
(\ref{gE11}):
\begin{align*}
[H^*(X;rat(M)]:=k_*[M]=k'_*f_*[M]=\sum_j(-1)^{j+n} \cdot
k'_*[H^{j+n}(f_*M)] \\ =\sum_j(-1)^j \cdot [H^*(Y;R^jf_*rat(M))] \in
K_0(mHs^p).
\end{align*}
Note that in the case when $M=\Q^H_X$ and $Y$ is compact, this
implies the claimed identity (\ref{E11}) of Prop. \ref{smooth}.

If we assume in addition that the local system $R^jf_*rat(M)$ is
constant with stalk the mixed Hodge structure $M^j$, then, by
``rigidity", when viewed as a mixed Hodge module it is isomorphic to
$k'^*M^j\simeq \Q^H_Y \otimes k'^*M^j$. Therefore, by the
$K_0(MHM(pt))$-linearity of $k'_*$, this implies
$$[H^*(Y;R^jf_*rat(M))]=[H^*(Y;rat(\Q^H_Y))] \cdot [M^j] \in K_0(mHs^p).$$
Next note that for a fixed $y \in Y$,
$$\sum_j(-1)^j \cdot [M^j]=[H^*(\{f=y\};rat(M))] \in K_0(mHs^p).$$
So, if all direct image sheaves $R^jf_*rat(M)$ $(j \in \Z)$ are
constant, then we obtain the multiplicative relation claimed in the
formula (\ref{gmult}):
$$[H^*(X;rat(M)]=[H^*(Y;\Q)] \cdot [H^*(\{f=y\};rat(M))] \in K_0(mHs^p).$$
By using $k'_!$ instead of $k'_*$ in the above arguments, we get a
similar multiplicative relation for the cohomology with compact
support.

\end{proof}

\subsection{The main theorem and immediate
consequences.}\label{MAIN}
We can now formulate the main technical result of this paper, and
explain, as an application, the proofs of Theorems \ref{formula1}
and \ref{formula2}

Let $Y$ be an irreducible complex algebraic variety endowed with a
complex algebraic Whitney stratification $\VV$ so that the
intersection cohomology complexes $$IC'_{\bar W}:=IC_{\bar W}[-{\rm
dim}(W)]$$ are $\VV$-constructible for all strata $W \in \VV$.
Denote by $S$ the top-dimensional stratum, so $S$ is Zariski open
and dense, and $V \leq S$ for all $V \in \VV$. Let us fix for each
$W \in \VV$ a point $w \in W$ with inclusion $i_w:\{w\}
\hookrightarrow Y$. Then
\begin{equation}\label{100} i_w^*[IC'^H_{\bar W}]=[i_w^*IC'^H_{\bar
W}]=[\Q^H_{pt}]\in K_0(MHM(w))=K_0(MHM(pt)),\end{equation} and
$i_w^*[IC'^H_{\bar V}] \neq [0] \in K_0(MHM(pt))$ only if $W \leq
V$. Moreover,  for any $j \in \Z$, we have
\begin{equation}\label{cone} \HC^j (i_w^*IC'_{\bar V}) \simeq
IH^j(c^{\circ}L_{W,V}),\end{equation} with $c^{\circ}L_{W,V}$ the
open cone on the link $L_{W,V}$ of $W$ in $\bar V$ for $W \leq V$
(cf. \cite{B}, p.30, Prop.4.2). So
$$i_w^*[IC'^H_{\bar V}]=[IH^*(c^{\circ}L_{W,V})] \in K_0(MHM(pt)),$$
with the mixed Hodge structures on the right hand side defined by
the isomorphism (\ref{cone}).

The main technical result of this section is the following
\begin{thm}\label{main}  For each stratum $V \in \VV \setminus
\{S\}$ define inductively
\begin{equation}\label{eq8}
\widehat{IC^H}(\bar V):=[IC'^H_{\bar V}] - \sum_{W < V}
\widehat{IC^H}(\bar W) \cdot i_w^* [IC'^H_{\bar V}] \in
K_0(D^bMHM(Y)).
\end{equation}
As the notation suggests, $\widehat{IC^H}(\bar V)$ depends only on
the complex algebraic variety $\bar V$ with its induced algebraic
Whitney stratification. Assume $[M] \in K_0(D^bMHM(Y))$ is an
element of the $K_0(MHM(pt))$-submodule $\langle [IC'^H_{\bar V}]
\rangle$ of $K_0(D^bMHM(Y))$ generated by the elements $[IC'^H_{\bar
V}]$, $V \in \VV$. Then the following equality holds in
$K_0(D^bMHM(Y))$:
\begin{equation}\label{mE}
[M]= [IC'^H_Y] \cdot i_s^*[M]+\sum_{V < S}  \widehat{IC^H}(\bar V)
\cdot \left( i_v^*[M] -i_s^*[M] \cdot i_v^*[IC'^H_Y] \right).
\end{equation}
\end{thm}

\begin{proof}
In order to prove formula (\ref{mE}), consider
\begin{equation}\label{coef}[M]=\sum_{V \in \VV} [IC'^H_{\bar V}] \cdot L(V),\end{equation}
for some $L(V) \in K_0(MHM(pt))$. The aim is to identify these
coefficients $L(V)$. Since $S$ is an open stratum, by applying
$i_s^*$ to (\ref{coef}) we obtain:
$$i_s^*[M]=L(S) \in K_0(MHM(s))=K_0(MHM(pt)).$$
Next fix a stratum $W\neq S$, and apply $i_w^*$ to (\ref{coef}).
Recall that $i_w^*[IC'^H_{\bar W}]=[\Q^H_{pt}]\in
K_0(MHM(w))=K_0(MHM(pt)),$ and $i_w^*[IC'^H_{\bar V}] \neq [0] \in
K_0(MHM(pt))$ only if $W \leq V$. We obtain
\begin{equation}\label{1} i_w^*[M]=L(W)+\sum_{W < V}
i_w^*[IC'^H_{\bar V}] \cdot L(V) \in
K_0(MHM(w))=K_0(MHM(pt)).\end{equation} Since $S$ is dense, we have
that $W < S$, so the stratum $S$ appears in the summation on the
right hand side of (\ref{1}). Therefore
\begin{equation}\label{2}
i_w^*[M]-i_w^*[IC'^H_{Y}] \cdot i_s^*[M]=L(W)+\sum_{W<V<S}
i_w^*[IC'^H_{\bar V}] \cdot L(V) \in K_0(MHM(w))=K_0(MHM(pt)).
\end{equation}
This implies that we can inductively calculate $L(V)$ in terms of
$$L'(W):=i_w^*[M]-i_w^*[IC'^H_{Y}] \cdot i_s^*[M].$$
Indeed, (\ref{2}) can be rewritten as \begin{equation}\label{3}
L'(W)=\sum_{W\leq V<S} i_w^*[IC'^H_{\bar V}] \cdot L(V) \in
K_0(MHM(pt)),
\end{equation}
and the matrix $A=(a_{W,V})$ with $a_{W,V}:=i_w^*[IC'^H_{\bar V}]
\in K_0(MHM(pt))$ for $W,V \in \VV \setminus \{S\}$ is
upper-triangular with respect to $\leq$, with ones on the diagonal.
So $A$ can be inverted. The non-zero coefficients of
$A^{-1}=(a'_{W,V})$ can inductively be calculated (e.g., see
\cite{Sta}, Prop. 3.6.2) by $a'_{V,V}=1$ and
\begin{equation}\label{4}
a'_{W,V}=-\sum_{W \leq T < V} a'_{W,T} \cdot a_{T,V}
\end{equation}
for $W<V$. Then equation (\ref{coef}) becomes
\begin{equation}\label{5}[M]=[IC'^H_Y] \cdot i_s^*[M] + \sum_{W < S}  [IC'^H_{\bar W}] \cdot L(W)=
[IC'^H_Y] \cdot i_s^*[M] + \sum_{W \leq V< S}  [IC'^H_{\bar W}]
\cdot a'_{W,V} \cdot L'(V).
\end{equation}
The result follows by the inductive identification (for $V<S$
fixed):
$$\sum_{W\leq V} [IC'^H_{\bar W}] \cdot
a'_{W,V}=[IC'^H_{\bar V}]-\sum_{W \leq T <V} [IC'^H_{\bar W}] \cdot
a'_{W,T} \cdot a_{T,V}=\widehat{IC}^H(\bar V).$$

\end{proof}

Before stating immediate consequences of the above theorem, let us
describe some cases when the technical hypothesis $[M] \in \langle
[IC'^H_{\bar V}] \rangle$  needed in the proof is satisfied for a
fixed $M \in D^b(MHM(Y))$.

\begin{example}\label{MainEx}\rm
\begin{enumerate}\item Assume that all sheaf complexes $IC'_{\bar V}$, $V \in \VV$, are
not only $\VV$-constructible, but satisfy the stronger property that
they are ``cohomologically $\VV$-constant", i.e., all cohomology
sheaves $\HC^j(IC'_{\bar V})|_W$ ($j \in \Z$) are constant for all
$V,W \in \VV$. Moreover, assume that either \begin{enumerate}\item
$rat(M)$ is also cohomologically $\VV$-constant, or
\item all perverse cohomology sheaves $rat(H^j(M))$ ($j \in \Z$) are
cohomologically $\VV$-constant, e.g., each $H^j(M)$ is a pure Hodge
module with $\HC^{-dim(V)}(rat(H^j(M))|_V)$ constant for all $V \in
\VV$.\end{enumerate} Then $[M] \in \langle [IC'^H_{\bar V}]
\rangle$. In particular, if all strata $V \in \VV$ are
simply-connected, then we have that $[M] \in \langle [IC'^H_{\bar
V}] \rangle$ for all $M \in D^bMHM(X)$ so that $rat(M)$ is
$\VV$-constructible.

\item \emph{Toric varieties.} Another interesting example comes from a toric variety $Y$
with its natural Whitney stratification $\VV$ by orbits, i.e., each
stratum is of the form $V \simeq (\C^*)^{dim (V)}$. In this case not
all the strata are simply connected, but nevertheless all
intersection complexes $IC'_{\bar V}$, $V \in \VV$ are
cohomologically $\VV$-constant, e.g., see \cite{BL}, Lemma 15.15. It
follows that for any $M \in D^b(MHM(Y))$ satisfying (a) or (b)
above, we have $[M] \in \langle [IC'^H_{\bar V}] \rangle$.

\begin{proof} It is clear that $(2)$ follows from $(1)$. Note that, by using the identity $$[M] =
\sum_i (-1)^i [H^i(M)] \in K_0(MHM(Y)),$$ it suffices to prove the
claim $(1)$ under the assumption $(a)$. Moreover, by using the
$t$-structure on $D^b(MHM(Y))$ that corresponds to the usual
$t$-structure on $D^b_c(Y)$ (and not to the perverse $t$-structure),
cf. [\cite{Sa1}, Remark 4.6(2)], we may assume that $rat(M)$ is a
sheaf, and by $(a)$ this is assumed to be cohomologically
$\VV$-constant (i.e., $rat(M)|_W$ underlies a constant variation of
mixed Hodge structures, for any $W \in \VV$).

Under these considerations, we can now prove the claim by induction
over the number of strata of the stratification. In the case when
$Y$ has only one stratum $U$, the variety $Y=U$ is non-singular,
connected, with $IC'^H_U = \Q^H_U$, and the claim follows trivially:
indeed, under our assumptions,  $M$ is just a constant variation of
mixed Hodge structures with stalk $F$, so $M \simeq \Q_Y^H \otimes
k^* F$, for $k:Y \to pt$ the constant map; therefore $[M]=[IC'^H_U]
\cdot [F]$. For the induction step, fix an open stratum $U \in \VV$
with open inclusion $j : U \hookrightarrow Y$ and closed complement
$i : Y' := Y \setminus U \to Y$. Then $U$ is a smooth, connected
variety as before, and $Y'$ is an algebraic variety with a smaller
number of strata, satisfying the same assumption on the intersection
cohomology sheaves of the strata. The distinguished triangle (cf.
\cite{Sa1}, (4.4.1))
$$j_!j^*M \to M \to i_*i^*M \to$$  implies that $$[M] = [j_!j^*M] +
[i_*i^*M] \in K_0(D^bMHM(Y )).$$ We can now apply the induction
hypothesis to $[i^*M]$. Moreover, by assumption, $rat(j^*M)$ is  a
constant sheaf with stalk the (graded polarizable) mixed Hodge
structure $F$, so that
$$j_![j^*M]=j_![j^*(IC'^H_{\bar U} \otimes k^*F)]=[IC'^H_{\bar
U}]\cdot [F] - i_*[i^*(IC'^H_{\bar U} \otimes k^*F)],$$ with $k: Y
\to pt$ the constant map. But
$$i^*(IC'^H_{\bar U} \otimes k^*F) \simeq i^*(IC'^H_{\bar U}) \otimes
i^*k^*F$$ has also a cohomologically $\VV$-constant underlying
complex, with respect to the induced stratification on $Y'$. So our
claim follows by induction.

\end{proof}

\end{enumerate}
\end{example}

\bigskip

In the following, we specialize to the relative context of a proper
algebraic map $f:X \to Y$ of complex algebraic varieties, with $Y$
irreducible, and indicate a proof of Theorem \ref{formula1} and of
Theorem \ref{formula2}. For a given $M \in D^bMHM(X)$, assume that
$Rf_*rat(M)$ is constructible with respect to the given complex
algebraic Whitney stratification $\VV$ of $Y$, with open dense
stratum $S$. By proper base change, we get
$$i_v^*f_*[M]=[H^*(\{f=v\}, rat(M))] \in K_0(MHM(pt)).$$
So under the assumption $f_*[M] \in \langle [IC'^H_{\bar V}]
\rangle$, Theorem \ref{main} yields the following identity in
$K_0(MHM(Y))$:
\begin{cor}\label{cormain}
\begin{eqnarray*}
f_*[M]&=& [IC'^H_Y] \cdot [H^*(F;rat(M))]  \\ &+& \sum_{V < S}
\widehat{IC^H}(\bar V) \cdot \left( [H^*(F_V;rat(M))]
-[H^*(F;rat(M))] \cdot [IH^*(c^{\circ}L_{V,Y})] \right),
\end{eqnarray*}
where $F$ is the (generic) fiber over the top-dimensional stratum
$S$, and $F_V$ is the fiber over a stratum $V \in \VV \setminus
\{S\}$.
\end{cor} Note that the corresponding classes $[H^*(F;rat(M))]$ and
$[H^*(F_V;rat(M))]$ may depend on the choice of fibers of $f$, but
the above formula holds for any such choice. If all strata $V \in
\VV$ are simply connected, then these classes are independent of the
choices made. By pushing the identity in Corollary \ref{cormain}
down to a point via $k'_*$, for $k':Y \to pt$ the constant map, and
using the fact that $k'_*$ is $K_0(MHM(pt))$-linear, an application
of the $\chi_y$-genus (which is a ring homomorphism) yields the
following:
\begin{prop}\label{gsmp} Under the above notations and assumptions, the following
identity holds in $\Z[y,y^{-1}]$:
\begin{eqnarray*}
&& \chi_y([H^*(X;rat(M)])= I\chi(Y) \cdot \chi_y([H^*(F;rat(M))]) \\
&+& \sum_{V < S} \widehat{I\chi}_y (\bar V) \cdot \left(
\chi_y([H^*(F_V;rat(M))]) -\chi_y([H^*(F;rat(M))]) \cdot
I\chi_y(c^{\circ}L_{V,Y}) \right).
\end{eqnarray*}
\end{prop}

\begin{remark}\rm
Note that Theorem \ref{formula1} follows from Proposition \ref{gsmp}
above if we take $M=\Q^H_X$. Similarly, for $X$ pure dimensional,
Theorem \ref{formula2} follows if we let $M=IC'^H_X$ (together with
the stalk identifications of Lemma \ref{id}).\end{remark}
\begin{remark}\rm Working with $k'_!$ in place of $k'_*$, a similar
argument yields the corresponding results for $\chi_y^c(\cdot)$ and
respectively $I\chi_y^c(\cdot)$. Moreover, an application of the
$E$-polynomials yields similar formulae for $IE(\cdot)$ and
$IE_c(\cdot)$.
\end{remark}

\section{Characteristic classes}\label{four}

In this section we construct a natural characteristic class
transformation, $MHT_y$, which for an algebraic variety $X$ yields a
twisted homology class $IT_y(X)$,  whose associated genus for $X$
compact is $I\chi_y(X)$. The main result of this section is a
formula for the proper push-forward of such a class, and is a direct
consequence of Theorem \ref{main} and Corollary \ref{cormain}. The
construction of $MHT_y$ follows closely ideas of a recent paper of
Brasselet, Sch\"urmann and Yokura (\cite{BSY}, Remark 5.3, 5.4, but
see also Totaro's paper \cite{To}, \S 7), and is based on Saito's
theory of mixed Hodge modules (cf. \S \ref{smhm}).

\subsection{Construction of the transformation $MHT_y$.}

For any $p \in \Z$, Saito constructed a functor of triangulated
categories
$$gr^F_pDR: D^bMHM(X) \to D^b_{coh}(X)$$ commuting with proper push-down.
Here $D^b_{coh}(X)$ is the bounded derived category of sheaves of
$\mathcal{O}_X$-modules with coherent cohomology sheaves. Moreover,
$gr^F_p DR(M)=0$ for almost all $p$ and $M \in D^bMHM(X)$ fixed. If
$\Q_X^H \in D^bMHM(X)$ is the constant Hodge module on $X$, and if
$X$ is smooth and pure dimensional then $gr^F_{-p} DR(\Q_X^H) \simeq
\Omega^p_X[-p]$.

The transformations $gr^F_pDR(M)$ are functors of triangulated
categories, so they induce functors on the level of Grothendieck
groups. Therefore, if $K_0(D^b_{coh}(X)) \simeq G_0(X)$  denotes the
Grothendieck group of coherent sheaves on $X$, we obtain the
following group homomorphism commuting with proper push-down:
\begin{equation}\label{grF}
MHC_*: K_0(MHM(X)) \to G_0(X) \otimes \Z[y, y^{-1}],
\end{equation}
$$[M] \mapsto \sum_p \left( \sum_i (-1)^{i} \HC^i ( gr^F_{-p} DR(M) ) \right) \cdot
(-y)^p.$$ Recall that $G_0$ is a covariant functor with respect to
the proper push-down $f_*$ defined as follows: if $f:X \to Y$ is an
algebraic map, then $f_* :G_0(X) \to G_0(Y)$ is given by
$f_*([\mathcal{F}]):=\sum_{i \geq 0} (-1)^i [R^if_* \mathcal{F}]$,
for $R^if_* \mathcal{F}$ the higher direct image sheaf of
$\mathcal{F}$.

Note also that by work of Yokura (cf. \cite{BSY} and references
therein) one can define the following  generalization of the
Baum-Fulton-MacPherson transformation for the Todd class:
\begin{equation}\label{td}
td_{(1+y)}:G_0(X) \otimes \Z[y, y^{-1}] \to H^{BM}_{2*}(X) \otimes
\Q[y, y^{-1}, (1+y)^{-1}],
\end{equation}
$$[\mathcal{F}] \mapsto \sum_{k \geq 0} td_k([\mathcal{F}]) \cdot (1+y)^{-k},$$
with $td_k$ the degree $k$ component of the Todd class
transformation $td_*$ of Baum-Fulton-MacPherson \cite{BFM}, which is
linearly extended over $\Z[y, y^{-1}]$. Since $td_*$ is degree
preserving, this new transformation also commutes with proper
push-down (which is defined by linear extension over $\Z[y,
y^{-1}]$).

\bigskip

We can now make the following definition (cf. \cite{BSY}, Remark
5.3):
\begin{defn}\rm
The transformation $MHT_y$ is  the composition of transformations:
\begin{equation}\label{IT}
MHT_y :=td_{(1+y)} \circ MHC_*: K_0(MHM(X)) \to H^{BM}_{2*}(X)
\otimes \Q[y,y^{-1},(1+y)^{-1}].
\end{equation}
By the above discussion, $MHT_y$ commutes with proper push-forward.
\end{defn}

\begin{example}\label{pt}\rm Let $\VB=((V_{\C},F), V_{\Q},K) \in MHM(pt)=mHs^p$. Then:
\begin{equation}\label{point}
MHT_y([\VB])=\sum_p td_0([gr^p_F V_{\C}]) \cdot (-y)^p = \sum_p
\text{dim}_{\C} (gr^p_F V_{\C}) \cdot (-y)^p = \chi_y([\VB]).
\end{equation}
Also, since over a point the twisted Todd transformation
$$td_{(1+y)}:K_0(pt)=\Z[y,y^{-1}] \to
\Q[y,y^{-1},(1+y)^{-1}]=H^{BM}_{2*}(pt;\Q)[y,y^{-1},(1+y)^{-1}]$$ is
just the identity transformation, we get that $MHC_*=MHT_y=\chi_y$
on $K_0(mHs^p)$.
\end{example}

\begin{defn}\rm
For a pure $n$-dimensional complex algebraic variety $X$, we define
\begin{equation}
IT_y(X):=MHT_y(IC'^H_X); \ \ \ IC_*(X):=MHC_*(IC'^H_X).
\end{equation}
\end{defn}

\begin{remark}\label{ty}\rm(\emph{Normalization})\newline
If $X$ is \emph{smooth and pure-dimensional}, then $IC'^H_X \simeq
\Q_X^H$ in $D^bMHM(X)$, and $$MHC_*(IC'^H_X)=\sum_p [\Omega_X^p]
\cdot y^p=:\lambda_y(T^*_X)$$  is the total $\lambda$-class of the
cotangent bundle of $X$. Therefore, by \cite{BSY}, Lemma 3.1,
$$IT_y(X)=td_{(1+y)}(\sum_p [\Omega_X^p] \cdot y^p)=T_y^*(TX) \cap
[X] =:T_y(X),$$ where $T_y^*(TX)$ is the modified Todd class that
appears in the generalized Hirzebruch-Riemann-Roch theorem, i.e.,
the cohomology class associated to the normalized power series
defined by $Q_y(\alpha):=\frac{\alpha (1+y)}{1-e^{-\alpha(1+y)}} -
\alpha y$. If $X$ is compact, the genus associated to $T_y^*(TX)$,
that is the degree of the zero-dimensional part of $T_y(X)$, is
exactly the Hirzebruch $\chi_y$-genus of Definition \ref{D1}.

\end{remark}

\begin{remark}\rm
It is conjectured in \cite{BSY} that for a complex projective
variety $X$, the homology class $IT_1(X)$ is exactly the $L$-class
$L_*(X)$ of Goresky-MacPherson. At least the equality of their
degree follows from Saito's work. In general, for any compact
complex algebraic variety $X$, one has that the degree of $IT_y(X)$
is exactly $I\chi_y(X)$, i.e.,
$$I\chi_y(X)=\int_X IT_y(X)$$
(e.g., see Corollary \ref{deg} below).
\end{remark}

\begin{remark}\rm
For a possibly singular algebraic variety $X$, the twisted homology
class $MHT_y(\Q_X^H)$ is the motivic Hirzebruch class $T_y(X)$
constructed in \cite{BSY,SY}, which for  a complete variety $X$ has
as associated genus the generalized Hirzebruch
$\chi_y$-characteristic from Definition \ref{D2} (cf. \cite{BSY}, \S
5). In particular, if $X$ is a rational homology manifold, then
$\Q_X^H \cong IC'^H_X$ in $D^b MHM(X)$ and $IT_y(X)=T_y(X)$.
Similarly, for any complex algebraic variety $X$, $MHC_*(\Q^H_X)$ is
the \emph{motivic Chern class} from \cite{BSY}.
\end{remark}

\subsection{Formula for proper push-forward}

We begin this section with the following simple
observation.\footnote{Finding numerical invariants of complex
varieties, more precisely Chern numbers that are invariant under
small resolutions, was Totaro's guiding principle in his paper
\cite{To}.} If $f:X \to Y$ is a proper algebraic map between
irreducible $n$-dimensional complex algebraic varieties so that $f$
is \emph{homologically small of degree $1$} in the sense of
\cite{GM2}, \S 6.2, then
$$f_*IT_y(X)=IT_y(Y).$$ Indeed, for such a map we have that $f_*IC_X
\simeq {^{p}\HC^0} (f_*IC_X) \simeq IC_Y$ in $D^b_c(Y)$ (\cite{GM2},
Theorem 6.2). Moreover, as $rat: MHM(Y) \to Perv(\Q_Y)$ is a
faithful functor, this isomorphism can be lifted to the level of
mixed Hodge modules (cf. \cite{Sa}, Theorem 1.12). Then, since
$MHT_y$ commutes with proper push-down and
$[IC'^H_X]=(-1)^n[IC^H_X]$ in $K_0(MHM(X))$, we obtain:
\begin{align*}f_*IT_y(X)=f_*MHT_y([IC'^H_X])=(-1)^nMHT_y(f_*[IC^H_X])
=(-1)^nMHT_y([IC^H_Y])\\ =MHT_y([IC'^H_Y])=IT_y(Y).\end{align*} In
particular, if $f:X \to Y$ is a \emph{small resolution}, that is a
resolution of singularities that is small in the sense of
\cite{GM2}, then:
\begin{equation}\label{small} IT_y(Y)=f_*T_y(X).\end{equation}
One might wish to take formula (\ref{small}) as a definition of the
class $IT_y(Y)$. Unfortunately, small resolutions do not always
exist.

\bigskip

The main result of this section is the following:
\begin{thm}\label{charfor}
Let $f :X \to Y$ be a proper  morphism of complex algebraic
varieties, with $Y$ irreducible. Let $\VV$ be the set of components
of strata of $Y$ in a stratification of $f$, with $S$ the
top-dimensional stratum (which is Zariski-open and dense in $Y$),
and assume $\pi_1(V)=0$ for all $V \in \VV$. For each $V \in \VV
\setminus \{S\}$, define inductively
$$\widehat{IT}_y(\bar V):= IT_y(\bar V)- \sum_{W < V}
\widehat{IT}_y(\bar W) \cdot I\chi_y(c^{\circ}L_{W,V}),$$ where
$c^{\circ}L_{W,V}$ denotes the open cone on the link of  $W$ in
$\bar{V}$, and all homology characteristic classes are regarded in
the Borel-Moore homology of the ambient variety $Y$ (with
coefficients in $\Q[y,y^{-1},(1+y)^{-1}]$). Then:
\begin{equation}\label{charformula}
f_*T_y(X)=IT_y(Y) \cdot \chi_y(F) + \sum_{V < S} \widehat{IT}_y(\bar
V) \cdot \left( \chi_y(F_V) - \chi_y(F) \cdot
I\chi_y(c^{\circ}L_{V,Y}) \right),
\end{equation}
where $F$ is the generic fiber of $f$, and $F_V$ denotes the fiber
over a stratum $V \in \VV \setminus \{S\}$.

If, moreover, $X$ is pure-dimensional, then:
\begin{equation}\label{charformula}
f_*IT_y(X)=IT_y(Y) \cdot I\chi_y(F) + \sum_{V < S}
\widehat{IT}_y(\bar V) \cdot \left( I\chi_y(f^{-1}(c^\circ L_{V,Y}))
- I\chi_y(F) \cdot I\chi_y(c^{\circ}L_{V,Y}) \right).
\end{equation}
\end{thm}

\begin{proof}
This follows directly by applying the transformation $MHT_y$ to the
identity of Corollary \ref{cormain} for $M=\Q^H_X$, and resp. for
$M=IC'^H_X$, and by using the fact that $MHT_y$ commutes with the
exterior product
$$K_0(MHM(Y)) \times K_0(MHM(pt)) \to K_0(MHM(Y\times \{pt\}))\simeq K_0(MHM(Y)).$$
More precisely, $MHT_y$ commutes with the first exterior product,
and with the last isomorphism induced by the proper push-down $p_*$
for the isomorphism $p:Y \times \{pt\} \overset{\simeq}{\to} Y$. If
$i$ is the inverse to $p$, and $k:Y \to pt$ is the constant map,
then  for $[M] \in K_0(MHM(Y))$ and $[M'] \in K_0(MHM(pt))$ we get
$$[M] \cdot [M']=[M \otimes k^*M']=[i^*(M \boxtimes M')]=[p_*(M \boxtimes M')].$$
Thus $$MHT_y([M] \cdot [M'])=p_*(MHT_y([M]) \boxtimes
MHT_y([M']))=MHT_y([M]) \cdot \chi_y([M']).$$

\end{proof}

\begin{cor}\label{deg} For any compact pure-dimensional complex algebraic variety $X$, the degree of $IT_y(X)$
is the intersection homology genus $I\chi_y(X)$, i.e., $$I\chi_y(X)=\int_X IT_y(X)$$
\end{cor}

\begin{proof} Apply Theorem \ref{charfor} to the constant map $f:X \to {point}$,
which is proper since $X$ is compact.

\end{proof}

\begin{remark}\rm Similar formulae can be obtained by applying the
transformation $MHC_*$ to the identity in Corollary \ref{cormain},
and even for a general mixed Hodge module $M$ so that $f_*[M] \in
\langle IC'^H_{\bar V}] \rangle$, e.g., if all strata $V \in \VV$
are simply connected.
\end{remark}

\begin{example}\label{bup2}\rm \emph{Smooth Blow-up}\newline
Let $Y$ be a smooth $n$-dimensional variety and $Z \overset{i}{
\hookrightarrow} Y$  a submanifold of pure codimension $r+1$. Let
$X$ be the blow-up of $Y$ along $Z$, and $f:X \to Y$ be the blow-up
map. Then $\VV:=\{Y \setminus Z,Z\}$ is a Whitney stratification of
$Y$ with $IC'_Y$, $IC'_Z$ and $Rf_*(\Q_X)$ all cohomologically
$\VV$-constant. Then as in Example \ref{bup}, we have that (compare
with \cite{BSY}, Example 3.3 (3)):
\begin{equation}\label{BSY2}
f_*T_y(X)=T_y(Y)+T_y(Z) \cdot \left( -y+\cdots +(-y)^r \right).
\end{equation}
In particular, for $y=0$, this yields the well-known formula
\footnote{This fact can be seen directly as follows: if $f:X \to Y$
is a blow-up map along a smooth center, then
$f_*([\mathcal{O}_X])=[\mathcal{O}_Y]$, for $f_*:G_0(X) \to G_0(Y)$
the proper push-forward on the Grothendieck groups of coherent
sheaves (e.g., see \cite{BW}, Lemma 2.2). The claim follows since
$td_*(X):=Td([\mathcal{O}_X])$, where $Td :G_0(X) \to
H_{2*}^{BM}(X;\Q)$ is the Baum-Fulton-MacPherson Todd transformation
\cite{BFM}, and the latter commutes with proper push-forward.}
(e.g., see \cite{BSY}, Example 3.3 (3)): $$f_*td_*(X)=td_*(Y).$$

\end{example}

\bigskip

As a special case, we also obtain the following generalization of
some well-known facts concerning multiplicative properties of
characteristic classes (e.g., see \cite{H}, \S 23.6, for a
discussion on Todd classes, or \cite{CS2} for a more general formula
for $L$-classes):
\begin{cor}\label{smoothchar} Let $f:X \to Y$ be a proper algebraic
map of complex algebraic varieties, with $Y$ smooth and connected,
so that all direct image sheaves $R^jf_*\Q_X$, or $R^jf_*IC'_X$ for
$X$ also pure-dimensional, are locally constant (e.g., $f$ is a
locally trivial topological fibration). Let $F$ be the general fiber
of $f$, and assume that $\pi_1(Y)$ acts trivially on the
(intersection-) cohomology of $F$ (e.g. $\pi_1(Y)=0$), i.e., all
these $R^jf_*\Q_X$ or $R^jf_*IC'_X$ are constant. Then:
\begin{equation} f_*T_y(X)=\chi_y(F) \cdot T_y(Y), \ \ \  f_*IT_y(X)=I\chi_y(F) \cdot T_y(Y).\end{equation}
\end{cor}

We conclude by pointing out that the assumption of trivial monodromy
is closely related, but different than the situation of ``algebraic
piecewise trivial" maps coming up in the motivic context (e.g., see
\cite{BSY}). For example, the first formula in Corollary
\ref{smoothchar} is true for a Zariski locally trivial fibration of
possibly singular complex algebraic varieties (see \cite{BSY},
Example 3.3).

\providecommand{\bysame}{\leavevmode\hbox
to3em{\hrulefill}\thinspace}

\end{document}